\documentclass{eurocg13}

\begin{document}
\title{New  representation  results for   planar
 graphs
\
}
\author{Farhad Shahrokhi\\
Department of Computer Science and Engineering,  UNT\\
P.O.Box 13886, Denton, TX 76203-3886, USA\
farhad@cs.unt.edu
}

\date{}
\maketitle
\date{} \maketitle

%%%%%%%%%%%%%%%%%%%%%%%%%%%%%%%%%%%%%%%%%%%%%%%%%%%%%%%%%%%%%%%%%%%%%%%%%%%%%%%% 

%%%%%%%%%%%%%%%%%%%%%%%%%%%%%%%%%%%%%%%%%%%%%%%%%%%%%%%%%%%%%%%%%%%%%%%%%%%%%%%% 
\begin{abstract}
A universal representation theorem is derived that shows any graph
is the intersection graph of  one chordal graph, a  number of
co-bipartite graphs, and one unit interval graph. Central to the 
the result is the notion of the  clique cover width which is 
a generalization of the bandwidth parameter. 
Specifically, we  show that  
any planar graph is
the intersection graph of one chordal graph, four  co-bipartite
graphs, and one  unit interval graph. 
Equivalently, 
any planar graph is
the intersection graph of a chordal graph and a graph that has 
{clique
cover width} of at most seven. 
We further
describe the extensions of the results to graphs drawn on  surfaces and  graphs
excluding a minor of crossing number of at most one.

\end{abstract}

\section{Introduction and Summary}
Graph theory, geometry, and topology
stem from the same roots. 
Representing  graphs as the 
intersection graphs
of geometric or combinatorial objects is highly desired   
in  certain branches of combinatorics, discrete and 
computational geometry, graph
drawing and information visualization, and the design of geographic information
systems (GIS). 
A suitable  
intersection model
not only provides a better understanding of the underlying graph, 
but it can also lead to computational advances.
A remarkable result in this area is  Koebe's 
(also Thurston's) 
theorem,  asserting
 that every planar graph is the touching graph of planar disks.
A similar result is due to Thomassen \cite{Th}   
 who showed  that every planar graph is the intersection graph of axis 
parallel boxes in $R^3$. 
Another noteworthy result is due to  Gavril \cite{gav}  
who proved that every chordal graph 
(a graph with no chordless cycles) is the intersection graph of
a collection of subtrees of a tree. 

Any (strict) partially  ordered set \cite{To1}  $(S,<)$ has a directed acyclic  graph
$\hat G$ associated with it in a natural way: $V(G)=S$, and $ab\in E(G)$
if and only if $a<b$. The {\it comparability graph}  associated with $(S,<)$
is  the undirected graph which is obtained by dropping the orientation on
edges of $\hat G$.
The complement of a comparability graph is an  {\it incomparability graph}.
Incomparability graphs are well studied due to their
rich structures and  are known to be the intersection  
graph of planar curves \cite{C1}. 
A interesting result in this area is due to Pach and 
T\"or\H ocsik 
\cite{PT} who  showed,  given a set of  straight line segments in the 
plane, there are four  incomparability graphs whose edge intersections
gives rise to the intersections of the segments.
Moreover, recent
work 
in combinatorial geometry 
has shown the connections 
between the intersection patterns of  arbitrary  planar curves
and properties of incomparability graphs  
\cite{FP1},   
\cite{FP2}. 

An  an interval graph is the intersection
graph of a set of intervals on the real line \cite{To2}. It is easily seen that
an interval graph is an incomparability graph. A unit interval graph is the
intersection graph of a set of unit intervals.

Throughout this paper, $G=(V(G), E(G))$ denotes  a connected undirected graph.
Let $d\ge 1$, be an integer,  and for $i=1,2,...,d$ let $H_i$ be 
a graph with $V(H_i)=V$, and  let $G$ be a graph with $V(G)=V$ and $E(G)=\cap_{i=1}^dE(G_i)$. 
Then  we say $G$ is the {\it intersection graph} of $H_1,H_2,...,H_d$, and write $G=\cap_{i=1}^tH_i$.
A clique cover $C$ in $G$ is a partition of $V(G)$ into cliques.  
 Throughout this paper,   we will  write
$C=\{C_0,C_1,...,C_t\}$ to indicate that $C$ is an ordered set of cliques.
For a clique cover $C=\{C_0,C_1,...,C_t\}$, in $G$, let
the {\it width}  of $C$, denoted by $W(C)$, denote
$\max\{|j-i||xy\in E(G), x\in C_i,y\in C_j,C_i,C_j\in C\}$.
The {\it clique cover width } of $G$ denoted by $CCW(G)$ is the smallest width
all ordered clique covers in $G$. Note that $CCW(G)\le BW(G)$, where 
$BW(G)$ denotes the bandwidth of $G$.
A co-bipartite graph is the complement of a bipartite graph. Clearly, 
any  co-bipartite graph is an incomparability graph.  
\subsection {Our Results}
We recently proved the following result \cite{Sh}.

\begin{theorem}\label{t2}
{\sl 
Let  $C$ be a clique cover in $G$ with $0<W(C)\le w, w\ge 1$.
Then, there  are
$\lceil\log({w})\rceil+1$ co-bipartite graphs $H_i, i=1,2,...,
\lceil\log({w})\rceil+1$,
and a  unit interval  graph $H_{\lceil\log({w})\rceil}+2$,
so that $G=\cap_{i=1}^{\lceil\log(w)\rceil +2}H_i$.
}
\end{theorem}

The main result in this paper  is Theorem \ref{t9}, which  asserts
any planar graph is the intersection graph of a chordal graph
and a graph whose clique cover width is bounded by seven. The application
of Theorem \ref{t2}, then, gives another version of the result as stated
in the abstract. Theorem \ref{t9} is obtained using the Universal
Representation Theorem, or Theorem \ref{t8}, which is interesting
on its own, and asserts that any graph is intersection graph of a chordal
graph and a graph  whose clique cover width is bounded. Nonetheless,
the upper bound on the clique cover width of the second graph 
depends on the properties of the tree decompositions  of the original graph.
Theorem \ref{t9} is further  extended to graphs drawn
on surfaces, and graphs excluding a minor with the crossing number of 
at most one.

%$\mu({\cal A}')=\mu({\cal A})$.

\section{Main Results}
\begin{defini}\label{d0}
{\sl 
A  {\it tree decomposition} \cite{RS}  
of a graph $G$ is a pair $(X,T)$
where $T$ is  a tree, and $X=\{X_i|i\in V(T)\}$ is a family of
subsets of $V(G)$, each called a bag, so that the following hold:
           
$\bullet$ $\cup_{i\in V(T)}X_i=V(G)$
            
$\bullet$ for any $uv\in E(G)$, there is an $i\in V(T)$ so that
$v\in X_i$ and $u\in X_i$.

$\bullet$ for any $i,j,k\in V(T)$, if $j$ is on the path from
$i$ to $k$ in $T$, then $X_i\cap X_k\subseteq X_j$.

}
\end{defini}

\begin{theorem}\label{t8}(Universal Representation Theorem)
{\sl
Let $G$ be a graph and let $L=\{L_1,L_2,...,L_k\}$ be  a
partition of vertices, so that for any  $xy\in E(G)$, either  $x,y\in L_i$
where $1\le i\le k$, or, 
$x\in L_i, y\in L_{i+1}$, where,   $1\le i\le k-1
$.  Let $(X,T)$ be a tree decomposition of $G$. Let
$t^*=\max_{i=1,2,...,k}\{|L_i\cap X_j||j\in V(T)\}$. (Thus, $t^*$ 
is the largest number of vertices in  any element of $L$  that appears in any  
bag of $T$). 
Then, there is a graph $G_1$ with $CCW(G_1)\le 2t^*-1$ and a chordal
graph $G_2$ so that  $G=G_1\cap G_2$.
}
\end{theorem}
{\bf Proof.}
For any $v\in V(G)$, let $X_v$ be the set of bags in $X$ that contains
vertex $v$,  and let $T_v$ be the subtree of $T$ on the vertex set $X_v$.   
Let $G_2$ be the intersection graph of these subtrees. Thus, for any
$v,w\in V(G)$, $vw\in E(G_2)$, if $X_v\cap X_w\ne\emptyset$. It is well known
that $G_2$ is chordal. See work of Gavril  \cite{gav}. 
Now let $\omega$ be the largest clique in $G_2$ among all cliques whose vertices
are entirely in $L_i$, for some $i=1,2,...,k$.
It follows from established properties on the tree decomposition that
all vertices in $\omega$ should appear in one bag $B$ in $X$. Consequently,
$|\omega|\le |B\cap L_i|\le t^*$. Next observe that   
for $i=1,2,...,k$, $G_2[L_i]$ is chordal and hence perfect,
and thus there must be at most $t^*$ disjoint independent  sets $L^j_i, j=1,2, ...,t^*$ whose union is $L_i$.
Now construct $G_1, V(G_1)=V(G)$, as follows: $E(G_1)=E(G)\cup E'$,
where $E'$ is obtained by placing an edge between any vertex pair in 
each independent set ${L_i}^j$ for $i=1,2,,...,k, j=1,2,...,t^*$.
Clearly, $G=G_1\cap G_2$. 
In addition,  for $i=1,2,...,k$, $G_1[L_i]$ is covered with at most $t^*$ disjoint
cliques, hence any  ordering of these cliques will give rise to a clique
cover $C$ of $G_1$ with $W(C)\le 2t^*-1$, since  any edge $e\in E(G_1)$
either has both ends in one previously prescribed clique in $G_1[L_i]$,
or must have end points in  two consecutive elements in $L$. $\Box$.

The following definitions are from \cite{bod}.

\begin{defini}\label{d3}
{\sl 
A {\it maximal spanning forest} 
of $G$ is a spanning forest $T$ that contains a spanning tree from
each component of $G$. Thus, when $G$ is connected, any spanning tree 
of $G$ is also a maximal spanning forest. Let $T$ be a maximal spanning
tree of $G$, and  
let $ab\in E(G)-E(T)$; The {\it detour} of  $ab$ in $T$ is  
 the unique $ab$ path in $T$. Let $e\in E(T)$. 
The {\it  edge remember number} of  $e$, 
denoted by $er(e,T,G)$, is the number of  edges in 
$E(G)-E(T)$ whose
detour contains $e$; Equivalently, $er(e,T,G)$ is the number of fundamental
cycles in $G$ relative to $T$, that contain $e$. Similarly, for $v\in V(G)$,  
the {\it  vertex remember number number} of $v$ denoted by
$vr(v,T,G)$, is the number of edges in $E(G)-E(T)$ whose detour, 
or the fundamental cycle associated with it,  contains $v$.
To remedy technical issues, for any $e\in E(G)-E(T)$, we define $er(e,T,G)=0$.
The edge  remember number  and vertex  
remember number of $G$ in $T$, 
denoted by $er(G,T)$ and $vr(G,T)$,   are the largest
remember numbers overall edges in $E(T)$ and vertices in $V(T)$, respectively.  

}
\end{defini}

\begin{defini}\label{d4}
{\sl
Let $T$ be a maximal spanning tree of $G$, and let  
$\hat T$ be a forest 
that is obtained by inserting vertices of degree two to the edges of $T$.
Thus, ${\hat T}=(V(T)\cup E(T), E({\hat T}))$.  
Now, for  any $v\in V({T})$ place $v$ in $X_v$, and for any $e=ab\in E(T)$ place
$a$ and $b$ in $X_{e}$. Next, for any $e=ab\in E(G)-E(T)$, take one of $a$ or
$b$, say $a$, and place it in $X_v$, for any $v$ which is on the unique
$ab-$detour in $T$; Similarly, place $a$ in $X_e$ for any edge $e$ which
is on the unique $ab-$detour in $T$. 
Finally, define,
${\hat X}=\{X_i|i\in V(T)\cup E(T)\}$.
}
\end{defini}

Bodlaender \cite{bod} showed the following.
\begin{theorem}\label{bod} 
{\sl Let  $T$ be a maximal spanning tree of  
$G$, and let $\hat T $ and $\hat X$  be as defined above. Then, 
$({\hat X}, {\hat T})$ is a tree decomposition of $G$ whose width 
is at most $\max\{vr(G,T),er(G,T)+1\}$.
}
\end{theorem}

In light of the above result, we will refer to $({\hat T}, {\hat X})$ (in  definition \ref{d4})   as 
a {\it tree decomposition
of $G$ relative to $T$}.  
Note that the construction in definition \ref{d4} would allow the 
same vertex to appear in $X_v$ or $X_e$ more than once, where each 
appearance is associated 
with  an end point of an edge $e\in E(G)-E(T)$, representing
a distinct fundamental cycle containing $v$, or, $e$. With that in
mind, we have , $|X_v|=vr(v,T,G)+1$ and $|X_e|=er(e,T,G)$+2.  
However, when viewing $|X_v|$ and $|X_e|$ as sets, the duplicate members 
would be removed, thereby, $=$ would become $\le$.

The following Lemma is extended from \cite{bod}. 
The notations and
claims are slightly perturbed to exhibit additional properties
of the construction of Bodlaender, that we will use later.  

\begin{lemma}\label{l10}
{\sl Let $G$ be a plane  graph, let $O$ be the set of all vertices
in the outer boundary of $G$, let  $H, V(H)=V(G)$ be  a graph
obtained by removing all edges in the outer boundary $G$.  Let 
$T'$ be a maximal spanning forest of $H$ and let 
$({\hat X}',{\hat T}')$ be 
a tree decomposition of $H$ relative to $T'$.

$(i)~~$ $T$ can be extended to a maximum  spanning forest $T$ of $G$
so that 
$vr(v,T,G)\le vr(v,T',H)+\Delta(G)$ and 
$er(e,T,G)\le er(e,T',H)+2$,   for all $v\in V(G)$ and $e\in E(T)$.

$(ii)~~$ $({\hat X}',{\hat T}')$ can be extended to a tree
decomposition $({\hat X}, {\hat T})$ of $G$ relative to $T$ so that 
$|X_v\cap O|\le |X'_v\cap O|+\Delta(G)$ and $|X_e\cup O|\le |X'_e\cup O|+2$
for all $v\in V(G)$ and $e\in E(T)$  
.\footnote{ In $(i)$ and $(ii)$ 
we follow the assumption that $er(e,T',G)=0$ and $X'_e=\emptyset$,  for $e\in E(T)-E(T')$.} 
}
\end{lemma}
{\bf Proof.}  
For $(i)$, let $K$ be graph with $V(K)=V(G)$ and $E(K)=E(T')\cup (E(G)-E(H))$, and
note that the external face of $K$ is the same as external face of $G$.
Extend $T'$ to a maximal spanning tree $T$ of $K$ by adding edges from
$E(G)-E(H)$. Note that for any $e=xy\in E(K)-E(T)$, $x$ and
$y$ must be on the boundary of $G$. Thus, the associated  $xy$ detour $p$ in $T$
plus $e$ must form the boundary of a non-external  face in $K$. Since any edge
in $T$ is  common  to at most 2 
non external faces,  and each vertex in $T$ is common  to 
at most $\Delta(G)$ many non-external faces, in $K$, it follows that  
for any $e\in E(T)$ and any $v\in V(G)$,  $er(e,T,K)\le 2$ and
$vr(v,K,T)\le \Delta(G)$. As $T$ is also a maximal spanning tree of $G$
and each fundamental cycle in $G$ is either a fundamental cycle of
$K$ relative to $T$, or a fundamental cycle of $H$ relative to $T'$,
we must have $er(e,G,T)\le er(e,T',H)+er(e,T,K)\le er(e,T'H)+2$, 
and $vr(v,G,T)\le er(e,T',H)+vr(v,T,K)\le vr(v,T',H)+\Delta(G)$.

$(ii)$ follows from $(i)$.  In particular, note that additional 2 
or $\Delta(G)$ fundamental edges  that 
contribute to $vr(v,G,T)$ and $er(e,G,T)$, respectively,  are those
edges in $E(G)-E(T)$ that have both end points in $O$. 
Now obtain a tree decomposition of $G$ relative to $T$, by extending 
each bag of  ${\hat T}'$, to a bag of ${\hat T}$ by the possible addition
f one end point of such a fundamental edge, as
described in definition \ref{d4}.  $\Box$ 

By a  plane graph we mean an embedding of a planar graph in the plane.
A plane  graph is  $1-$outer planar, if it is outer planar.
For $k\ge 2$, a plane  graph $G$ is $k-$outer planar, if after  removal of
all vertices (and edges incident to these vertices) in the external face of 
$G$, a  $k-1$outer planar graph is obtained.

\begin{theorem}\label{t9}
{\sl
Let $G$ be a planar graph,
then, there is a graph $G_1$ with $CCW(G_1)\le 7$ and a chordal
graph $G_2$ so that  $G=G_1\cap G_2$.
}
\end{theorem}

{\bf Proof.} 
Assume $G$ is $k-$outer planar. Thus,  there are graphs 
$G=G_1,G_{2},...,G_k$
so that for $i=1,2,...,k$,
$G_i$ is $(k-i+1)-$outer planar, and $G_{i+1}$ is obtained by removing the vertices 
in 
the outer face of $G_{i}$. 
For $i=1,2,...,k$, let
$O_i$ denote the set of  vertices on the outer 
face of $G_i$. Note that for $i=1,2,...,k$,
one can replace any vertex $v$ of degree $d\ge 4$ in the outer face of 
 $O_i$ 
by a path $p_v$  of   
$d-2$ vertices of degree 3, so that $G$ is transformed to another 
$k-$outer planar graph $G'$.  Specifically, for $i=1,2,...,k$, let 
$O'_i$ denote the 
set of vertices corresponding to $O_i$, after   this transformation. 
Note that $G'$ is $k-$outer planar and has maximum degree 3,
let $G'_1=G'$, and  for $i=2,...,k+1$,  
let $G'_{i}$ denote the graph that is 
obtained after removing all edges in  the outer face of $G'_{i-1}$, 
and note that $G'_{i}$ is $(k-i+1)-$outer planar and of maximum degree
3. Note that $G'_{k+1}$
is acyclic and let  $T_{k+1}=G'_{k+1}$.  
Clearly, 
$vr(v,T_{k+1},T_{k+1})=0,er(e,T_{k+1},T_{k+1})=0$,  
for any $x\in V(G)$, and any $e\in E(T_{k+1})$.  
Thus, for the
tree decomposition $({\hat X_{k+1}}, {\hat T_{k+1}})$ of $G_{k+1}$ relative 
to $T_{k+1}$, 
and bags $X_v, X_e$, 
$v\in V(G), e=ab\in E(T_{k+1})$, we have  
$|X_v|=1$ (since $X_v=\{v\}$), and $|X_e|=2$ (since $X_e=\{a,b\}$),  
respectively.
Next, for $j=k,k-1,...1$, let  $T_j$ and $({\hat X_j}, {\hat T_j})$ be a maximal spanning forest and  a tree decomposition of $G'_j$ relative to $T_{j}$, that are obtained by the application
of Part $(i)$ and Part $(ii)$ of  Lemma  \ref{l10}, to $T_{j+1}$ and 
$({\hat X_{j+1}}, {\hat T_{j+1}})$, respectively. 
Thus, $({\hat X_{1}}, {\hat T_{1}})$ is a tree decomposition of $G'$.
Then, one can show (by induction) that  
for any $j,i=k,k-1,...,1$, and any $X^j_v, X^j_e\in {\hat X_j}$ with  
$v\in V(G), e\in E(T_j)$ 

$|X^j_v\cap O'_i|=|X^{j-1}_v\cap O'_i| \hbox{~if~~} i\ne j, 
\hbox{whereas,~} |X^j_v\cap O'_i|\le 1+\Delta(G'_j)\le 1+3=4  \hbox{~if~} i=j,$

and

$|X^j_e\cap O'_i|=|X^{j-1}_e\cap O'_i| \hbox{~if~} i\ne j, 
\hbox{~whereas,~} |X^j_e\cap O'_i|\le 2+2=4  \hbox{~if~} i=j.$

Hence, for $i=1,2,...,k$, 
and $X^1_v, X^1_e\in {\hat X}_1$ with
$v\in V(G)$ and  $e\in E(T_1)$, we have, 
$|X^1_v\cap O'_i|\le 4$ and  $|X^1_e\cap O'_i|\le 4$. 
Next,  for any $v\in V(G)$,  contract
all the vertices in $p_v$ to $v$, thereby, for $i=1,2,...,k$  contracting
 $O'_i$ to $O_i$. 
For 
any bag $X^1_{t}\in {\hat X}_1$ with  $t\in V(G)\cup E(T_1)$, let 
$Y_t=(X^1-p_v)\cup \{v\}$.
Now let ${Y}=\{Y_t|t\in {V(G)\cup E(T_1)}\}$. Since $G$ is a minor of $G'$, 
it follows that $(Y,T^1)$ is a tree decomposition of $G$  with the property
that
for any $Y_t\in {Y}$ with  $t\in V(T_1)\cup E(T_1)$, and
any  $i=1,2,...,k$, 
we have 
$|Y_t\cap  O_i|\le 4$. 
Now the result follows from Theorem \ref{t8}, by taking 
$L=\{O_1,O_2,...,O_k\}$. 
$\Box$.

Combining Theorems \ref{t2} and \ref{t9} we obtain the following.
\begin{theorem}\label{t10}
{\sl
Let $G$ be a planar graph,
then, there are co-bipartite graphs $G_1,G_2,G_3,G_4$, a unit interval graph $G_5$, 
and a chordal graph $G_6$ so that  $G=\cap_{i=1}^6G_i$. 
}
\end{theorem}
\subsection{Extensions}
The result for  planar graphs give rise to the following.

\begin{theorem}\label{t11}
{\sl Let $G$ be a graph of genus $g$. Then, there is an integer   
$c=O(\log(g))$, co-bipartite graphs $G_i,i=1,2,...,c$, a unit interval
graph $G_{c+1}$, and a chordal graph $G_{c+2}$  so that $G=\cap_{i=1}^{c+2}G_i$.
}
\end{theorem}
{\bf Proof Sketch.}
One can show the claim by induction on $g$, where Theorems \ref{t9} and  \ref{t10}
establish the base of the induction.
$\Box$

\begin{theorem}\label{t12}
{\sl
Let $G$ be a graph that does not have as a minor, a graph $H$ 
whose crossing number is at most one. 
Then there is an integer $c=O(\log(C_H))$,  
co-bipartite graphs $G_i,i=1,2,...,c$,  a unit interval
graph $G_{c+1}$, and a chordal graph $G_{c+2}$ so that $G=\cap_{i=1}^{c+2}G_i$.
Here, $C_H={20}^{2{(2|V(H_0)|+4|E(H_0)|)}^5}$.
}
\end{theorem}
{\bf Proof Sketch. } It is known that any graph that does not have a minor $H$ of 
crossing number of at most one,  can be obtained by taking the  clique sum
of a finite set of graphs, where each graph is either planar, or has a 
tree width of at 
most $C_H$ \cite{RS}. 
So $G=H_1\bigoplus H_2...\bigoplus H_k$, where $\bigoplus$ stands for 
the clique sum operation, and 
for $i=1,2,...,k$, each $H_i$ is either planar, or has a tree width of 
at most $C_H$. 
We prove the  claim by induction on $k$. 
When  $k=1$ the result follows from Theorems \ref{t9}, \ref{t8}, \ref{t2}, 
and the definition of $C_H$. Now assume that the claim is true for $k=t-1$,
let $k=t\ge 2$, and set $F=H_1\bigoplus H_2...\bigoplus H_{t-1}$. 
Then, $G=F\bigoplus H_k$. By induction,  $F=F_1\cap F_2$, where  
$F_2$ is chordal and $CCW(F_1)\le C_H$. Moreover,  since $H_{t-1}$
is either planar, or has a tree width of at most $C_H$, by Theorem \ref{t9} 
we have  $G_{t-1}=F_3\cap F_4$, where $CCW(F_3)\le 2C_H$  and $F_4$ is chordal. 
Now let $G_1=F_1\bigoplus F_3$, and $G_2=F_2\bigoplus F_4$, then, it
is easy to verify that $G_2$ is chordal. To finish the proof, one can verify
using properties of the clique cover width that,  $CCW(G_2)\le 2C_H$.
Now the claim follows from Theorem \ref{t2}.
$\Box$

\section{Computational   Aspects}
All constructions provided here can be done in polynomial time, with
the exception of Theorem  \ref{t12}.  
%The details will be provided,
%if the paper is accepted for presentation.

In \cite{Sh} we have shown that if   $G$ is the 
intersection graph of a chordal graph and a graph whose clique cover width is
bounded by a constant, 
then  $G$ can be separated with a splitting ratio of $1/3-2/3$, 
for a variety of measures, where the measure associated with the separator
is ``small". Consequently, the planar separator theorem \cite{LT}
and its extensions follow from the representation  results in this paper.

We  highly suspect that the computation of the clique cover width
is an NP$-$hard problem,  due to its connection with the bandwidth problem.


\begin{thebibliography}{99}

\small






\bibitem{Bo}
Bodlaender H.L,  A Tourist Guide through Treewidth. Acta Cybern. 11, 1993,  1-22.

\bibitem{bod}
 Bodlaender H.,  A partial k-arboretum of graphs with bounded treewidth. Theoretical Computer Science, 209, 1998, 1-45.



\bibitem{C1}  Golumbic M.,  Rotem D.,  Urrutia J.,  Comparability graphs and intersection graphs, Discrete Mathematics 43 (1), 1983, 37-6.




\bibitem{Ch}
Chan T.,
Polynomial-time approximation schemes for packing and piercing fat objects
, Journal of Algorithms,  46(2), 2003, 178 - 189. 


\bibitem{FP2} Fox J. and  Pach J.,
String graphs and incomparability graphs,
Advances in Mathematics, 2012, 1381-1401.



\bibitem{FP1}
Fox J.,   Pach J.,
A separator theorem for string graphs and its applications,
 Combinatorics, Probability and Computing 19, 2010, 371-390.


%\bibitem{FPT} Fox J, Pach J,  T\'oth C.D, Intersection patterns of curves,
%J. London Math. Soc. (2008)


\bibitem{gav}  Gavril, F.,  The intersection graphs of subtrees in trees 
are exactly the chordal graphs, Journal of Combinatorial Theory, Series 
B 16, 1974, 47-56. 

\bibitem{LT} Lipton R. J., Tarjan R.E. , A separator theorem for planar graphs, SIAM Journal on Applied Mathematics 36, 1979,  177-189



\bibitem{PT}  
Pach J., 
T\"or\H ocsik J., 
Some geometric applications of Dilworth's theorem,
{\it Disc. Comput. Geometry}, 21,  1994, 1-7.

\bibitem{RS}
 Robertson N.,  Seymour, P.  D.  Graph minors III: Planar tree-width, 
Journal of Combinatorial Theory, Series B 36 (1), 1984, 49-64.

\bibitem{Sh}
 Shahrokhi F., in preparation. 




\bibitem{Th} Thomassen, C., Interval representations of 
planar graphs, 
Journal of Combinatorial Theory, Series B 40, 1986, 9-20.

\bibitem{To2} Trotter W.T.,   New perspectives on interval orders and interval
graphs,  in Surveys in Combinatorics, Cambridge Univ. Press, 1977, 237-286.

\bibitem{To1}  Trotter, W.T., Combinatorics and partially ordered sets: Dimension theory, Johns Hopkins series in the mathematical sciences, The Johns Hopkins University Press, 1992.





\end{thebibliography}
\end{document}